\documentclass[11pt,twoside, amscd]{article}

\usepackage{amssymb,amsmath,amsthm}

\newcommand{\barefootnote}[1] {%
  \begingroup
    \renewcommand{\thefootnote}{}
    \footnotetext{#1}
    \renewcommand{\thefootnote}{\arabic{footnote}}
  \endgroup
}

\renewcommand{\title}[1] {%
  \begingroup
    \begin{center}
      \vspace{0.4in}
      \bf\huge
      \addtolength{\baselineskip}{5mm}
      #1
    \end{center}
  \endgroup
}

\newcommand{\url}[1] {%
  \barefootnote{%
    {\small e-print archive: }
    {\texttt http://xxx.lanl.gov/#1}
  }
}

\renewcommand{\author}[1] {%
  \begingroup
    \begin{center}
      \vspace{0.4in}
      \bf
      #1
      \vspace{0.2in}
    \end{center}
  \endgroup
}

\newcommand{\address}[1] {%
  \begingroup
    \begin{center}
      #1
    \end{center}
 \endgroup
}

\newcommand{\addressemail}[1] {%
  \begingroup
    \begin{center}
      \vskip-\baselineskip
      #1
    \end{center}
  \endgroup
}

\newcounter{secondpage}
\newcounter{gfirstpage}
\newcounter{glastpage}
\newcommand{\copyrightnotice}[4]{
  \leftline{\copyright~ #1 International Press}
  \leftline{Adv. Theor. Math. Phys. {\bf #2} (#1) #3--#4}
  \setcounter{gfirstpage}{#3}
  \setcounter{glastpage}{#4}
}
\begin{document}


\pagestyle{myheadings}
\thispagestyle{empty}


\newcommand{\cX}{{\cal X}}
\newcommand{\cY}{{\cal Y}}
\newcommand{\cR}{{\cal R}}
\newcommand{\cS}{{\cal S}}
\newcommand{\cT}{{\cal T}}
\newcommand{\cU}{{\cal U}}
\newcommand{\cV}{{\cal V}}
\newcommand{\cB}{{\cal B}}
\newcommand{\cC}{{\cal C}}
\newcommand{\cE}{{\cal E}}
\newcommand{\cO}{{\cal O}}
\newcommand{\cL}{{\cal L}}
\newcommand{\cM}{{\cal M}}
\newcommand{\cA}{{\cal A}}
\newcommand{\cF}{{\cal F}}
\newcommand{\cH}{{\cal H}}
\newcommand{\cf}{{\cal f}}
\newcommand{\cGamma}{{\cal\Gamma}}
\newcommand{\MKB}{{$M_{(K+B)^{2}}$}}
\newcommand{\cMC}{{{\cM\cC}} }
\newcommand{\bMC}{{{\bold M \bold C}}}
\newcommand{\cMCS}{{$\cM\cC(\cS)$}}
\newcommand{\XBL}{{$(X,B,L)$}}
\newcommand{\cXBL}{{${(\cX,\cB,\cL)}$}}
\newcommand{\cXBLo}{{$(\cX_{0},\cB_{0},\cL_{0})$}}
\newcommand{\cXBLs}{{$(\cX_{s},\cB_{s},\cL_{s})$}}
\newcommand{\Span}{\operatorname{Span}}
\newcommand{\Chow}{\operatorname{Chow}}
\newcommand{\codim}{\operatorname{codim}}
\newcommand{\supp}{\operatorname{supp}}
\newcommand{\Hilb}{\operatorname{Hilb}}
\newcommand{\Pic}{\operatorname{Pic}}
\newcommand{\red}{\operatorname{red}}
\newcommand{\Aut}{\operatorname{Aut}}
\newcommand{\PGL}{\operatorname{PGL}}
\newcommand{\Spec}{\operatorname{Spec}}
\newcommand{\Sym}{\operatorname{Sym}}
\newcommand{\Sing}{\operatorname{Sing}}
\newcommand{\id}{\operatorname{id}}
\newcommand{\diff}{\operatorname{d}}
\newcommand{\Rad}{\operatorname{Rad}}
\newcommand{\Ker}{\operatorname{Ker}}

\newtheorem{thm}{Theorem}[section]
\newtheorem{question}[thm]{Question}
\newtheorem{assumption}[thm]{Assumption}
\newtheorem{lem}[thm]{Lemma}
\newtheorem{cor}[thm]{Corollary}
\newtheorem{cors}[thm]{Corollaries}
\newtheorem{prop}[thm]{Proposition}
\newtheorem{crit}[thm]{Criterion}
\newtheorem{conj}[thm]{Conjecture}
\newtheorem{defn}[thm]{Definition}
\newtheorem{condition}[thm]{Condition}
\newtheorem{say}[thm]{}
\newtheorem{exmp}[thm]{Example}
\newtheorem{prob}[thm]{Problem}
\newtheorem{alg}[thm]{Algorithm}
\newtheorem{rem}[thm]{Remark}           
\newtheorem{note}[thm]{Note}            
\newtheorem{summ}[thm]{Summary}         
\newtheorem{ack}{Acknowledgments}       
\newtheorem{notation}{Notation}         
\newtheorem{notationnum}[thm]{Notation}
\newtheorem{claim}[thm]{Claim}
\newtheorem{case}{Case}
\newtheorem{subcase}{Subcase}
\newtheorem{step}{Step}
\newtheorem{approach}{Approach}
\newtheorem{principle}{Principle}
\newtheorem{fact}{Fact}
\newtheorem{subsay}{}
\setcounter{section}{0}

\def\proof{\smallskip{\it Proof. }}
\def\endproof{\hfill\qed}
\def\qed {\nobreak$\quad$\lower 1pt\vbox{\hrule
\hbox to 8pt{\vrule height 8pt\hfil\vrule height 8pt}
\hrule}\ifmmode\relax\else\par\medbreak \fi}
\def\exact#1#2#3{0\rightarrow{#1}\rightarrow{#2}\rightarrow{#3}\rightarrow0}



 \copyrightnotice{2002}{6}{557}{574}   
\setcounter{page}{557}
\title{HyperK\"ahler Manifolds  and Birational Transformations}

\url{math.AG/0111089}       

\author{Y. Hu and S.-T. Yau}

\address{The University of Arizona, Tucson, AZ \addressemail{yhu@math.arizona.edu}}
 \address{Harvard University, Cambridge, MA  \addressemail{yau@math.harvard.edu}}


\markboth{HYPERK\"AHLER MANIFOLDS  AND BIRATIONAL GEOMETRY }{Y. HU AND S.-T. YAU}

\section{Introduction}
\label{sec:introduction}

The minimal model program has greatly enhanced
our knowledge on birational geometry of varieties of dimension 3 and higher.
About the same time, the last two decades have also witnessed increasing
interests in HyperK\"ahler manifolds, a particular class of Calabi-Yau manifolds.
One interest in this area, which we hope
to treat in the future, is to investigate the behavior of the SYZ mirror conjecture
(\cite{SYZ})
under birational maps. The present paper focuses
on the birational geometry of projective symplectic varieties
and  attempts to reconcile the two approaches
altogether in {\it arbitrary} dimensions.

Let $X$ be a smooth projective
symplectic variety of dimension $2n$
with a holomorphically symplectic form $\omega$. Hence $X$ admits a HyperK\"ahler structure
by a theorem of Yau. For an arbitrary variety $V$, let  $V_{reg}$ denote the
non-singular part of $V$. Also, by a small birational contraction,
we will mean a non-divisorial birational contraction.

A striking and simple relation between minimal model program and symplectic
geometry is that a small Mori contraction necessarily coincides
(generically) with the null foliation of the symplectic form along the exceptional locus.

\bf Theorem A.
\it Let $\pi: X \rightarrow Z$ be a small birational contraction,
$B$ an arbitrary irreducible component of the degenerated locus
and $F$  a generic  fiber  of 
  $\pi: B  \rightarrow S:= \pi (B)$. Then
$T_b F = (T_b B)^\perp$ for  generic  $b \in F_{reg} \cap B_{reg}$. In particular,
\begin{enumerate}
\item the inclusion $j: B \hookrightarrow X$ is a coisotropic embedding;
\item the null foliation of $\omega|_B$ through
generic points of $B$ coincides with generic fibers of  $\pi: B  \rightarrow S$.
\end{enumerate}
\rm

This bears an interesting consequence.

\bf Theorem B.
\it Let $\pi: X \rightarrow Z$ be a small birational contraction and $B$ the exceptional locus.
Assume that $\pi: B \rightarrow S$ ($:= \pi (B)$) is a smooth fibration.
 Then  $\pi: B \rightarrow S $
is a ${\Bbb P}^r$-bundle over $S$ where $r = \codim B \ge 2$.
\rm

This shows, for example,  that the dimensional condition
that Mukai posed to define his elementary transformation
is necessary.

Specializing to the case when $S$ is a point, we obtain

\bf Theorem C.
\it Let $B$ be a smooth subvariety of a projective symplectic variety $X$ of dimension $2n$.
Assume that $B$ can be contracted to a point. Then
$B$ is isomorphic to  ${\Bbb P}^n$.
\rm

In particular, $B$ must be a Lagrangian subvariety. Indeed, this is a general phenomenon.

\bf Theorem D.
\it Let $B$ be an arbitrary subvariety of a projective symplectic variety $X$.
Assume that $B$ can be contracted to a point. Then
$B$ is Lagrangian.
\rm

More generally, we prove (consult \S \ref{section:semismall} for the explanation of the relevant
technical terms)

\bf Proposition E.
\it A small  birational contraction
$\pi: X \rightarrow Z$ is always IH-semismall.
\rm

Indeed, we expect that a stratum is IH-relevant to $\pi$ if and only if it carries the reduced
symplectic form $\omega_{red}$. In particular, $\pi$ is strictly IH-semismall if and only if
$Z$ is symplectically stratified using reduced symplectic forms.
In fact, we believe that  all small symplectic birational contractions are strictly IH-semismall.
(Conjecture \ref{strictlysemismall}).

The next relation between Mori theory and
symplectic geometry is that a Mori contraction necessarily produces a large family of Hamiltonian
flows around the exceptional locus, which makes the contraction
behaving like a holomorphic moment map locally around the exceptional locus.
This was observed jointly with Dan Burns.
The ``moment map'' depends on the singularities. It becomes more complicated
as the singularity becomes worse. In any case,
the fiber of the contraction over a singular point is preserved by
the corresponding Hamiltonian flows. Under some technical assumption
we were abel to see that  the generic fiber
of $B \rightarrow S$ is a finite union of almost homogeneous spaces.
In fact, we anticipate that every fiber should be a finite union of homogeneous spaces and
a generic fiber is  the projective space whose dimension coincides with
the codimension of the exceptional locus.

However, in general, it is
possible that these (almost) homogeneous spaces are not
generated by the above Hamiltonian flows, some other vector fields are needed.
The Halmitonian approach should be useful in the absence of projective structure.

Finally we conjecture that any two smooth symplectic varieties that are birational
to each other are related by a sequence of Mukai elementary transformations after removing
closed subvarieties of codimension greater than 2 (Conjecture \ref{codim2met}).

We end the introduction with the following question.

Let $K_0$ be the Grothendick ring of  varieties over ${\Bbb C}$: this is the abelian group
generated by isomorphism classes of complex varieties, subject to the relation $[X-S] = [X] - [S]$,
where $S$ is closed in $X$.
All the explicitly known symplectic birational maps are of symmetric nature. We wonder whether
the following is true:
Let $X$ and $Y$ be smooth symplectic varieties. Then
$[X]=[Y]$ if and only if $X$ and $Y$ are birational.

 In dimension 4, it follows from  \cite{BHL} and \cite{WW} that
$[X]=[Y]$ if  $X$ and $Y$ are birational. The same is true for smooth Calabi-Yau 3-folds.

The current paper is sequel to \cite{BHL} but does not depend on
it. We mention that we were aware of the work of
Cho-Miyaoka-Shepherd-Barron part of which are parallel to ours;
Only until August 2002, we were informed the existence of the
paper of Cho-Miyaoka-Shepherd-Barron when the first author was
visiting Taiwan. After that we learnt from the referee that that
paper has now appeared as \cite{Miyaoka3}. After our paper was
circulated, the paper of Wierzba and Wisniewski (\cite{WW})
appeared in ArXiv (math.AG/0201028). Our paper was circulated in
November 2001 (math.AG/0111089), and the results were reported
earlier at Harvard differential geometry seminars in Fall 1999.
There are also papers \cite{Namikawa}, \cite{wierzba} and
\cite{keladin} that deal with related topics with some overlaps.
Similar results to Theorem D and Proposition E were also obtained
in \cite{Namikawa}.

\medskip
{\sl Acknowledgments.} YH is very indebted to  Professor J\'anos
Koll\'ar for kindly answering every of his questions. He thanks
Professor Dan Burns for sharing many of his ideas back several
years ago. He also acknowledges the financial support (summer
1999) from MPI in Bonn while part of this paper was being prepared
and reported at the MPI Oberseminar. Last but not the least YH
thanks Harvard university and the second author for hospitality
and financial support (1999-2000) which has made this work
possible. The authors are partially supported by NSF  grants DMS
0124600 and  DMS 9803347.

\section{Deformation of rational curves}

Rational curves,   occurring in the exceptional loci of our birational
transformations, are studied in this section to prepare for what follows.
The results are, however,  of  independent interest.

A rational curve $C$ in a projective variety $W$ is called incompressible if
it can not be deformed to a reducible or non-reduced curve.

\begin{lem} (\cite{Miyaoka})
\label{rationkey}
Let $C$ be an incompressible rational curve
in a projective algebraic variety  $W$ of dimension $n$. Then  $C$
 moves in a family of dimension at most $2n-2$.
\end{lem}
\proof
Replacing $W$ by its normalization if necessary, we may assume that $W$ is normal. Then
the lemma is equivalent to (1.4.4) of \cite{Miyaoka}. (One needs to observe that up to
(1.4.4), the smooth assumption on the ambient variety is not used.)
\endproof

If there is a morphism $f: {\Bbb P}^n \rightarrow W$ such that it
is an isomorphism or a normalization, then $2n-2$ is minimal and achieved.

\begin{thm} (\cite{Miyaoka})
\label{miyaoka}
Let $W$ be a smooth projective  variety  of dimension $n$.
 Assume that $W$ contains a rational curve and every rational curve of $W$
moves in a family of dimension at least $2n -2$. Then $W \cong {\Bbb P}^n$.
\end{thm}
\proof We borrow some proofs from \cite{Miyaoka}.

If $\Pic (W) \ge 2$, there is a non-trivial extremal contraction
$$\gamma: W \rightarrow W'$$ such that $\dim W' > 0$. Let $S \subset \Chow (W)$
be an irreducible component of the set of the contracted extremal rational curves of
minimal degree and $F_S$ the incident variety over $S$.
Then for a general point $x \in  p_2 (F_S)$, $F_S (x) $, the universal
family of curves in $F_S$ that pass the point $x$,  has its image contained in a fiber of $\gamma$.
That is,
$$ \dim p_2 (F_S (x)) \le n-1.$$
By (1.4.4) of \cite{Miyaoka},
$$\dim S \le \dim p_2 (F_S) + \dim p_2 (F_S (x)) -2 \le n + n-1 -2 = 2n -3.$$
This is a contradiction.

Hence $\Pic (W) = 1$. $W$ is therefore  Fano. Then the result of \cite{Miyaoka} applies\footnote{
In the summerof 1999. J\'anos Koll\'ar informed YH of the paper \cite{Miyaoka} where
the theorem was proved for $W$ being smooth and Fano. We observed this theorem much later
(although independently). Our assumption on the number of the family is by Riemann-Roch
a consequence of their numerical assumption $-K_W \cdot C \ge 2n-2$. We makes no mention of $-K_W$
in the statement. Note that a nice rigorous proof of Cho-Miyaoka-Shepherd-Barron's theorem
was recently provided by Kebekus \cite{kebekus}}.
\endproof

In fact, we expect

\begin{conj}
Let $W$ be an arbitrary  projective variety of dimension $n$.
Assume that $W$
contains a rational curve and every rational curve of $W$
moves in a family of dimension at least $2n -2$. Then
 there is a birational dominating morphism  ${\Bbb P}^n \rightarrow W$
such that it is either an isomorphism or normalization.
\end{conj}

We learned that this would follow from a recent work by Cho-Miyaoka-Shepherd-Barron.
When $n=2$, it can be proved as follows.

\begin{prop} (Koll\'ar.)
\label{kollarp2}  Let $S$ be a normal surface, proper over $\Bbb C$.
Then $S$ satisfies exactly one of the following:
\begin{enumerate}
\item  Every morphism $f: {\Bbb P}^1 \rightarrow S$ is constant;
\item  There is a morphism $f: {\Bbb P}^1 \rightarrow S$ such that
       $f$ is rigid.
\item  $S \cong {\Bbb P}^2$;
\item  $S \cong {\Bbb P}^1 \times {\Bbb P}^1$, or $S$ is isomorphic to
       a minimal ruled surface over a curve of positive genus
       or a minimal ruled surface with a negative section
       contracted.
\end{enumerate}
\end{prop}

\proof  If we have either of (1) or (2), we are done.
Otherwise, there is a morphism $f: {\Bbb P}^1 \rightarrow S$ deforms
in a 1-parameter family, thus $S$ is uniruled.

Let $p: \bar S \rightarrow S$ be the minimal desingularization
with the exceptional curve $E$. $\bar{S}$ is also uniruled,
hence there is an extremal ray $R$. There are 3 possibilities for $R$.

\begin{enumerate}
\item  $\bar{S} \cong {\Bbb P}^2$, thus also $S \cong {\Bbb P}^2$ (which implies (3)).
\item  $\bar{S}$ is a minimal ruled surface (which implies (4)).
\item  $R$ is spanned by a (-1)-curve $C_0$ in $\bar S$.
\end{enumerate}

But (3) is impossible, because the image of $C_0$ in $S$ would have
been rigid. The proof goes as follows. Assume the contrary that
$f_0: {\Bbb  P}^1 \cong C_0 \subset \bar S \rightarrow S$ is not rigid
and let $f_t: {\Bbb P}^1 \rightarrow S$ be a 1-parameter deformation.
For general $t$, $f_t$ lifts to a family of morphisms
$\bar{f}_t : {\Bbb P}^1 \rightarrow \bar S$. As $t \to 0$,
the curves $\bar{f}_t ({\Bbb P}^1)$ degenerate and we obtain a cycle
$$\lim_{t \to 0} \bar{f}_t ({\Bbb P}^1) = C_0 + F$$
where $\text{ Supp} F \subset \text{ Supp} E$.
$\bar{S}$ is the minimal resolution, thus $K_{\bar{S}} \cdot F \ge 0$.
Therefore,
$$K_{\bar{S}} \cdot \bar{f}_t ({\Bbb P}^1) \ge K_{\bar{S}} \cdot C_0 = -1.$$
On the other hand, for a general $t$ the morphism $\bar{f}_t$ is free,
thus
$$K_{\bar{S}} \cdot \bar{f}_t ({\Bbb P}^1) \le -2$$
by II.3.13.1, \cite{Ko}. This contradication shows that
$f_0$ is rigid.
\endproof

\begin{lem}
\label{contractionlemma}
Let $B$ be a compact subvariety in a smooth complex variety. Assume that
$B$ is contractible and $C$ a compact curve in $B$ whose image is a point. Then
$C$ can not be moved out of $B$.
\end{lem}
\proof (J\'anos Koll\'ar.) Assume otherwise, there would be compact curves outside $B$
which are arbitrarily close to the curve $C$. Let $\pi: X \rightarrow Z$ be the contraction.
Then images of these curves form a family of curves shrinking to the point $\pi(C)$ in $Z$.
This is impossible because it would imply that a Stein neighborhood of  $\pi(C)$ contains
compact curves.
\endproof

\section{Symplectic varieties}

An arbitrary smooth variety is called
(holomorphically) symplectic if it carries a closed holomorphic two form
$\omega$ such that it is non-degenerated.

Our objective is to understand relations between any two symplectic models within the
same birational class.  Hence we only consider {\it small} birational contractions.
For simplicity, generic fiber of a birational contraction is sometimes assumed to be irreducible.

One of our first observations is

\begin{prop}\footnote{Thanks are due to Conan Leung whose comments
lead to this proposition.}
Let  $\Phi: X ---> X'$ be  a birational map between two
smooth symplectic varieties such that it is isomorphic over
open subsets $U \subset X$ and $U' \subset X'$, then for any
 symplectic form $\omega$ over $X$, there is a unique symplectic form
$\omega'$ over $X'$  such that $\omega|_U = \Phi^* \omega'|_{U'}$.
\end{prop}
\proof
Since $\Phi$ is necessarily isomorphic in codimension one, every closed form  over $U'$
extends to a closed form over $X'$. Let $\omega'$ be the extention to $X'$ of
$(\Phi^{-1})^* (\omega|_{U})$. The question is whether this extension is necessarily non-degenerated. Let $\dim X = 2n$. Then $(\omega')^n$ is a section of a line bundle over $X'$
so that its zero locus is a divisor if not empty. Since $\codim (X' \setminus U') \ge 2$,
$\omega'$ has to be non-degenerated. The uniqueness is clear.
\endproof

Mukai discovered a class of well-behaved birational maps
between symplectic varieties.

\begin{defn}
\label{mukai}
Let $P$ be a closed subvariety of $X$. Assume that
\begin{enumerate}
\item $P$ is ${\Bbb P}^r$-bundle over a symplectic manifold $S$;
\item $r = \codim P$.
\end{enumerate}
Then $X$ can be blown up along $P$ and the resulting variety can be
contracted along another ruling, yielding a symplectic variety $X'$.
Moreover, $X' \setminus U = P^*$ is the dual bundle of $P$
and $X'$ carries a symplectic form $\omega'$ which agrees with $\omega$ over $U$.
This  process is known as a Mukai's elementary transformation. Sometimes we abbreviate it
as MET.
\end{defn}

The following is a very important and useful property of a symplectic variety.

\begin{lem}  (Ran)
\label{ran} Let $W$ be an arbitrary smooth symplectic variety of dimension $2n$. Then
every proper rational curve in $W$ deforms in a family of dimension  at least $2n-2$.
\end{lem}
\proof
When the curve is smooth, it was proved by Ran.
His proof can be extended to singular case by  using graphs
(see e.g.,  Lemma 2.3 of \cite{BHL}).

Ran's original proof uses Block's semiregularity map.
There is a second proof without using the Block semiregularity map.
By Fujiki \cite{fujiki}, there is a (general) deformation of $X$ so that
no curves (indeed no subvarieties) survive in the nearby members. Let $\hat{X}$ be
the total space of such a family and $C$ a rational curve in the central fiber $X$.
Then $K_{\hat{X}} \cdot C = K_X \cdot C = 0$. Hence by Riemann-Roch, we have
$$\dim \Hilb_{[C]} \ge 2n + 1 -3 = 2n -2.$$
\endproof

Combining this lemma and Lemma \ref{rationkey}, we obtain that

\begin{prop} Let $\Phi: X ---> X'$ be a birational map between two smooth projective
symplectic varieties and $E$ the exceptional locus. Then
$$ n \le \dim E \le 2n -2.$$
\end{prop}

\begin{exmp}
Let $S$ be an elliptic fibration K3 surface with a $A_2$-type singlar fiber, that is,
a singular fiber is a union of two smooth rational curves $C_i, i=1,2$, crossing at
two distinct points P and Q. Let $X=S^{[2]}$. Then $X$ contains two $P^2s$, $C_i^{[2]}, i=1,2$,
meeting at the point $P +Q \in X=S^{[2]}$.  We can flop $C_i$ to get $X_i$, $i=1,2$, respectively.
Then the exceptional locus of
$$X_1 ---> X_2$$
is a union of $P^2$ and the Hirzebruch surface ${\Bbb F}_1$ intersecting
along the $(-1)$-section $C$ of  ${\Bbb F}_1$.

Note that $C$ deforms in a 1-dimensional family in  ${\Bbb F}_1$, hence must deform out
out of ${\Bbb F}_1$ by Lemma \ref{ran}. The strictly transform of this family in $X$
under the birational map $X_1 ---> X$ is a deformation family of a nodal curve $\ell_1 + \ell_2$
passing the point $P+Q$ with a branch (line) $\ell_i$ in each of $C_i^{[2]}, i=1,2$.  This nodal
curve  $\ell_1 + \ell_2$ can be deformed out of $C_1^{[2]} \cup C_2^{[2]}$.

\end{exmp}

\section{Mori contraction and null foliation}

Let $j: P \hookrightarrow X$ be a subvariety of $X$. Write $\omega|_P := j^* \omega$.
When $P$ is singular, by the form $\omega|_P$ we mean the restriction to the nonsingular part
$P_{reg}$ of $P$.

\begin{lem}
\label{lemma:noforms}
Let $\pi: X \rightarrow Z$ be a small birational contraction and $F$ a fiber of $\pi$.
Then $H^{0, q} (F) = 0$ for all $q > 0$. In particular, $F$ does not admit nontrivial
global holomorphic forms.
\end{lem}
\proof
Since Z has (at worst) canonical singularities,  it has rational singularities
 by Corollary 5.24 of \cite{km}. That is, $R^q \pi_* \cO_X = 0$ for all $q > 0$.
This implies that  $H^q (F, \cO_F)=0$ for all $q > 0$. Let $F'$ be a smooth resolution of $F$.
Then $H^q (F', \cO_{F'})=0$, hence $H^0 (F', \Omega^q_{F'})=0$ for all $q> 0$.
This implies that $F$ does not admit global holomorphic
forms because otherwise the pullback of a global form from $F$ to $F'$ yields a
contradiction to the fact that $H^0 (F', \Omega^q_{F'})=0$ for all $q> 0$.
\endproof

Unless otherwise stated, we will use $E$  to denote the degenerated
locus of the map $\pi$ and $B$ to denote an arbitrary irreducible component of $E$.
For a smooth point $b$ of $B$, $(T_b B)^\perp$  denotes the orthogonal complement
of the tangent space $T_bB$ with respect to the form $\omega_b$.

\begin{thm}
\label{omega=0} Let $\pi: X \rightarrow Z$ be a small birational contraction,
$B$ an arbitrary irreducible component of the degenerated locus
and $F$  a generic  fiber  of 
  $\pi: B  \rightarrow S:= \pi (B)$. Then
$T_b F = (T_b B)^\perp$ for a generic  $b \in F_{reg} \cap B_{reg}$ such that
$s=\pi (b)$ is a smooth point in $S$. In particular,
\begin{enumerate}
\item the inclusion $j: B \hookrightarrow X$ is (generically) a coisotropic embedding;
\item the null foliation of $\omega|_B$ through
generic points of $B$ coincides with generic fibers of  $\pi: B  \rightarrow S$.
\end{enumerate}
\end{thm}
\proof
First we have $\omega |_F = 0$ because $H^{0,2}(F)=0$.

Let $s=\pi (b)$ for some $b \in F$ as stated in the statement of
the theorem.
For any $v \in T_s S$, let $v'$ be any lifting to $T_b B$.
Define $$\iota_v \omega (w) := \omega (w, v')$$
 for any $w \in T_b F$.
This is independent of the choice of the lifting because two choices
are differed by an element of $T_b F$ and $\omega|_F = 0$.
Thus  $\iota_v \omega$ is well-defined holomorphic one form on $F$.
By Lemma \ref{lemma:noforms}, $\iota_v \omega$ is identical to zero.
This shows that $T_b F \subset (T_bB)^\perp $.
In particular, $\dim F \le \codim B$ because $\omega$ is non-degenerate.

On the other hand, if $\dim F < \codim B$, that is, $\dim B < 2n - \dim F$,
then $\dim S = \dim B - \dim F < 2n - 2\dim F$, or $\dim S + 2 \dim F -2 < 2n -2$.
By e.g. \cite{km}, $F$ contains rational curves.  Take a rational curve $C$ in $F$
that satisfies the statement of Lemma \ref{rationkey}, then in a neighborhood of
$F$ in $X$, $C$ moves
in a family of dimension $\le$
$$\dim S + 2 \dim F -2 < 2n -2,$$
contradicting to Lemma \ref{ran}.

Hence $\dim F = \codim B$. This implies that $T_bF = (T_bB)^\perp$.

The rest of the statements follows immediately.
\endproof

\begin{cor}
\label{genericfibers} Let $F$ be a generic fiber of $\pi: B \rightarrow S$. Then
 $$\codim B = \dim F.$$
\end{cor}

In particular,  the assumption (2) in Definition \ref{mukai}
for Mukai's elementary transformations is not
only sufficient but also {\text necessary}.

Note also that Corollary  \ref{genericfibers} implies that

\begin{cor} Let the notation and assumptions be as before. Then
$$\dim F \ge 2$$
for every fiber $F$ of the projection $\pi: B \rightarrow S$.
\end{cor}

\begin{cor}
\label{p2}
Assume that the generic fiber $F$ of the birational contraction is an
irreducible  surface. Then
the normalization of $F$ is isomorphic to ${\Bbb P}^2$.
\end{cor}
\proof
This follows by combining the fact that $F$ contains a rational curve
and every rational curve in $F$ has to move in an (at least) 2-dimensional family and
Proposition \ref{kollarp2}.
\endproof


We believe more is true but are not able to give a proof.

\begin{conj}
Assume that the generic fiber $F$ is a surface. Then $F$ is
normal and hence isomorphic to ${\Bbb P}^2$.
\end{conj}

\begin{thm}
\label{pr} Assume that the generic fiber $F$ is smooth.
Then $F$ is isomorphic to ${\Bbb P}^r$ where $r = \dim F = \codim B$.
\end{thm}
\proof
Let $r = \dim F$. Then,
$F$ contains rational curve and every rational curve moves in a family of dimension
at least $2r -2$ to satisfy the fact that every rational curve in $B$ moves
in a family of dimension
at least $2n -2$. By Theorem \ref{miyaoka}, $F \cong {\Bbb P}^r$,
and by Theorem \ref{genericfibers}, $r= \codim B$.
\endproof

\begin{thm}
\label{reduction}
Let $S^0$ be the (largest) smooth open subset of $S$ such that
the restriction of $\pi$ to $$B^0=\pi^{-1}(S^0) \rightarrow S^0$$
is a fibration. Then
 $S^0$ carries a symplectic form $\omega_{red}$
induced from $\omega$ such that $\pi^* \omega_{red} = \omega|_{B^0}$.
\end{thm}
\proof
First note that the fibration $B^0=\pi^{-1}(S^0) \rightarrow S^0$ coincides with the null
foliation of the form $\omega|_{B^0}$ by Theorem \ref{omega=0}.
We can introduce a local coordinate $z_1, \ldots, z_m$ in $B^0$
about any smooth point of $B^0$
such that the leaves of the foliation are (locally) given by
$$z_1=\text{const},  \cdots,  z_k = \text{const}$$
and so the tangent space to the foliation is spanned by ${\partial / \partial z_{k+1}},
\ldots, {\partial / \partial z_{r}}$ at each point in the coordinate neighborhood. If we now
write
$$\omega|_B = \sum a_{ij} d z_i \wedge d z_j,$$
where $a_{ij}$ are functions. The condition that
$$\iota_{\partial / \partial z_{k+1}} \omega|_B
= \cdots = \iota_{\partial/ \partial z_{r}} \omega|_B =0$$ implies that $a_{ij} =0$ if
$i$ or $j$ is greater than $k$. The condition $d \omega|_B =0$ implies that $a_{ij}$
are functions of $z_1, \ldots, z_k$ only. This exhibits
the local neighborhood as a product and  the 2-form splits according to the product
as $\omega|_B + 0$. This  implies that
$\omega$ descends to a non-degenerate closed two form $\omega_{red}$ on $S^0$
such that  $\pi^* \omega_{red} = \omega|_{B^0}$.
\endproof

We will call $\omega_{red}$ the reduced symplectic form over $S^0$.

The foliation is in general only locally trivial in \'etale or analytic topology.
In the case that $X$ is the moduli spaces $M_H ({\bf v})$ of stable coherent sheaves
over a fixed K3 surface with a polarization $H$ and a Mukai vector ${\bf v}$,
Markman has constructed such examples (\cite{markman}).

It is useful to record an immediate corollary.

\begin{cor}
\label{fibration}
If $\pi: B \rightarrow S$ is a smooth fibration, then it coincides
with the integrable null foliation of $\omega|_B$. In particular, $\pi: B \rightarrow S$
is a ${\Bbb P}^r$-bundle where $r=\codim B$ and $S$ is a
smooth projective symplectic variety with the induced symplectic form $\omega_{red}$.
\end{cor}

\begin{rem}
\label{disjoint}
Let $\pi: X \rightarrow Z$ be a small birational contraction.
Let $B_i$ ($i \in I$) be all the irreducible
component of the degenerate locus such that
$B_i \rightarrow S_i:= \pi (S_i)$ is a smooth fibration for all $i$.
Then we believe $B_i$ ($i \in I$) must be mutually disjoint.
In particular, the Mukai elementary transformations can be performed along
these components in any order. This observation is, however, not needed later.
\end{rem}

\section{Lagrangian subvarieties and exceptional sets}

\begin{lem}
\label{pt-lag}
Let $B$ be a (possibly singular, reducible) subvariety of $X$.
Assume that $B$ is contractible. Then
$B$ is Lagrangian if and only if it is contractible to a point.
\end{lem}
\proof Assume that $B$ is Lagrangian.
Let $\pi: X \rightarrow Z$ be the contraction. If $\dim \pi(B) >1$, choose a small
neighborhood $U$ in $Z$ of a smooth point $s$ in $\pi(B)$. Let $g$ be a holomorphic function over $U$
such that $d g_s \ne 0$ along $\pi (B)$. Let $f$ be the pullback of $g$, defined over a neighborhood of
the fiber $F=\pi^{-1}(s)$. Then the Hamiltonian vector field $\xi$ defined by
$\iota_{\zeta_f} \omega = d f$ has the property that $\omega (\zeta, T_bB) = d f (T_bB)
=d g (d \pi (T_bB)) \ne 0$ for a generic point $b \in B$.
Since $B$ is Lagrangian, $\zeta_f$ is not tangent to $B$, hence deforms a generic
fiber $F$ around $s \in S$ outside $B$,
a contradiction. 

Conversely, assume that $B$ is contracted to a point. Then  $\omega|_B =0$
by Lemma \ref{lemma:noforms}.
That $\dim B > n$ is impossible because $\omega|_B = 0$ and $\omega$ is non-degenerated.
If $\dim B=m < n$, by Lemma \ref{rationkey},
 some rational curve will only move in a family of dimension bounded above
by $2m -2$, contradicting with Lemma \ref{ran}.
\endproof

\begin{rem}
Note that the proof of necessary part is purely analytic and hence the corresponding
statement is valid in the analytic category as well. In the projective category
and assuming that the contraction $\pi: X \to Z$ is small, this lemma is a direct
corollary of Theorem \ref{omega=0}.
\end{rem}

\begin{thm}
\label{smooth-pt-pn}
Let $P$ be a smooth closed  subvariety of $X$. Then
$P$ can be contracted to a point if and only if $P$ is isomorphic to ${\Bbb P}^n$.
\end{thm}
\proof If $P$ is isomorphic to ${\Bbb P}^n$, being rational,
therefore  Lagrangian, its normal bundle is
isomorphic to $T^* {\Bbb P}^n$, hence negative. Thus $P$ is  contractible to a point.
Conversely, if $P$ can be contracted to a point,
by Lemma \ref{pt-lag} it is a Lagrangian. Combining Lemma \ref{ran} and Theorem  \ref{miyaoka},
we conclude that $P \cong {\Bbb P}^n$.
\endproof

In general, we expect

\begin{conj}
Let $X^{2n} \rightarrow Z$ be a birational contraction such that the degenerate locus $B$ is
contracted to isolated  points. Then $B$ is a disjoint union of
 smooth subvarieties that are isomorphic
to ${\Bbb P}^n$.
\end{conj}

See \cite{BHL} for an affirmative answer when $n=2$ under the normality assumption. We were
recently informed that this was proved for $n=2$ in \cite{WW}. Also \cite{WW} indicated that
the general case has been  treated by Cho-Miyaoka-Shepherd-Barron.

\section{semismall and symplectic contractions}
\label{section:semismall}

Let $\pi: U \rightarrow V$ be a proper birational morphism.
Assume that $V= \bigcup_\alpha V_\alpha$ is
a stratification of $V$ by smooth strata
such that $\pi: U \rightarrow V$ is weakly stratified with respect
to this stratification. That is, $\pi^{-1}(V_\alpha) \to V_\alpha$ is a fibration
for all $\alpha$. Such a stratification always exsits.

\begin{defn}
\label{defn:small}
 $\pi: U \rightarrow V$ is called IH-semismall if for any $v \in V_\alpha$ we have
$$ 2\dim \pi^{-1}(v) \le \codim V_\alpha.$$
The stratum $V_\alpha$ is called IH-relevant to the map $\pi$
if the equality holds.
$\pi$ is  strictly semismall
if all strata $V_\alpha$ are  relevant to $\pi$.
It is known that the definition  is independent of the choice of a stratification.
\end{defn}

$\pi$ is IH-small if none of the strata is relevant to $\pi$ except the open stratum.
But by Corollary \ref{genericfibers} IH-small maps do not  exist for contractions
of symplectic varieties.

\begin{prop}
\label{semismall} Let $\pi: X \rightarrow Z$ be a contraction and
$B$  an arbitrary irreducible component of the exceptional locus.
Assume that $S:=  \pi (B) = \bigcup_\alpha S_\alpha$ is a stratification such that
$\pi: B \rightarrow  \bigcup_\alpha S_\alpha$ is weakly stratified.
Let $F_\alpha = \bigcup_i F_{\alpha, i}$ be the union of the irreducible components of
a fiber $F_\alpha$ over $S_\alpha$. $B_\alpha = \pi^{-1}(S_\alpha) \cap B$.
Then we have $$T_b F_{\alpha, i} \subset (T_b B_\alpha)^\perp$$ for generic $b \in (F_{\alpha, i})_{reg}
\cap (B_\alpha)_{reg}$. In particular, $\pi: X \rightarrow Z$ is always IH-semismall.
\end{prop}
\proof
Consider the fibration
 $\pi: B_\alpha \rightarrow S_\alpha$. Let $F_{\alpha, i}$ be an irreducible
component of the fiber $F_\alpha$. Since $H^{0, 2} (F_{\alpha, i}) = 0$
(Lemma \ref{lemma:noforms}), we have $\omega|_{F_{\alpha, i}}=0$.
Then the essentially same proof  as in the proof of
Theorem \ref{omega=0} will yield  that
$$T_b F_{\alpha, i} \subset (T_b B_\alpha)^\perp$$ for generic point $b \in F_{\alpha, i}$. This implies
that
$$\dim F_{\alpha, i} \le \dim  (T_b B_\alpha)^\perp $$
for all $i$. But
 $$ \dim  (T_b B_\alpha)^\perp= \dim X - \dim B_\alpha = \dim X - \dim F_\alpha - \dim S_\alpha$$
 Hence we have
$$2 \dim F_\alpha \le \codim S_\alpha$$
for all $\alpha$. This means that $\pi: X \rightarrow Z$ is IH-semismall.

Here is an easy observation

\begin{prop}
Every contraction $\pi: X \rightarrow Z$ can be made to be
strictly semismall  by removing subvarieties of
codimension greater than 2, while
still keeping the properness of the map.
\end{prop}
\proof
First note that for each irreducible component $B$ of the exceptional locus
 we have $\codim B \ge 2$.
By Theorem \ref{genericfibers}, over the generic part of each irreducible component
$B$, $\pi$ is strictly
semismall. Hence
 $\pi$  is always strictly semismall
after removing  subvarieties of codimension greater than 2 and still keeping the properness of
the map.
\endproof

\begin{defn}
\label{defn:symaps}
 Let $\pi: X \rightarrow Z$ be a  birational contraction and $\pi: B \rightarrow S$
be as before (cf. Theorem \ref{omega=0}). Let $S = \cup S_\alpha$ be a stratification
such that $B \to S$ is weakly stratified.
If every $S_\alpha$ admits a symplectic form $\omega_{red}$
 induced from $\omega|_{B_\alpha}$ (cf. Theorem \ref{reduction}),
 we will say that $Z$ admits a (reduced) symplectic stratification
and that the contraction $\pi$ is symplectically stratified.
\end{defn}

The two definitions \ref{defn:small} and \ref{defn:symaps}
are conjecturally to be related as follows.

Let $\pi: X \rightarrow Z$ be a contraction and
$B$  an arbitrary irreducible component of the exceptional locus.
Assume that $S:=  \pi (B) = \bigcup_\alpha S_\alpha$ is a stratification such that
$\pi: B \rightarrow  \bigcup_\alpha S_\alpha$ is weakly stratified.
Let $F_\alpha = \bigcup_i F_{\alpha, i}$ be the union of the irreducible components of
a fiber $F_\alpha$ over $S_\alpha$. $B_\alpha = \pi^{-1}(S_\alpha) \cap B$.
Then by Proposition \ref{semismall} we have,
$$T_b F_{\alpha, i} \subset (T_b B_\alpha)^\perp$$ for generic $b \in (F_{\alpha, i})_{reg}
\cap (B_\alpha)_{reg}$.
This inlusion becomes equality for  a generic fiber of $\pi: B  \rightarrow S:= \pi (B)$
(Theorem \ref{genericfibers}).
We expect the following is  true:
{\it the equality holds if and only if $S_\alpha$ carries a reduced symplectic form.
In particular, $\pi$ is strictly semismall if and only if
it can be  symplectically stratified.}

In Markman's paper (\cite{markman}), all the contractions he constructed
are indeed strictly semismall. This re-enforces our conjecture

\begin{conj}
\label{strictlysemismall}
Every birational contraction from a smooth projective symplectic variety
is necessarily strictly semismall. Equivalently, it can always be symplectically stratified.
\end{conj}


\section{Further Observation and speculations}

In this section we will try to explain that Hamiltonian flows might be useful in the study
of birational geometry of symplectic varieties in the purely complex context.

Consider and fix a generic fiber $F := F_\alpha$  over a generic (smooth) point $s$ of $S_\alpha$.  Let $B_\alpha := \pi^{-1}(S_\alpha)$.
Take a sufficiently small
open neighborhood $U_s$ of $s$ in $Z$ such that it is embedded in ${\Bbb C}^N$
where $N$ is the minimal  embedding dimension of $Z$ at $s$. Let $U_F = \pi^{-1} (U_s)$.
The composition of the maps
$$ U_F @> {\pi} >> U_s \hookrightarrow {\Bbb C}^N$$
is still denoted by $\pi$ (confusion does not seem likely).

Let $P$ be any subvariety of $U_F$.
By saying that a vector field $\xi$ on $U_F$ is tangent to $P$, we mean
$\xi$ is tangent to the nonsingular part of $P$. For a point $x \in U_F$, $\xi (x)$ denotes
the vector of the field at $x$.

\begin{prop}
\label{biggroupinfini}
For every  function $g$ about $s$ in ${\Bbb C}^N$,
let $f = \pi^* g$. Then the Hamiltonian vector field
$\zeta_f$ defined by $\iota_{\zeta_f} \omega = d f$ is tangent to $B_\beta \cap U_F$
for all $\beta$.
\end{prop}
\proof This is because otherwise
we would be able to deform (an irreducible component of)
$B_\beta \cap U_F$ along the flow of $\zeta_f$,
which is impossible by the weak stratification.
\endproof

\begin{prop}
\label{smallgroupinfini}
For every function $g$ about $s$ in ${\Bbb C}^N$
such that the differential $d g$ vanishes along $U_s\cap S$ at $s$, let $f=\pi^* g$.
Assume that $S_\alpha$ carries an induced symplectic form.
Then the Hamiltonian vector field
$\zeta_f$ defined by $\iota_{\zeta_f} \omega = d f$ is tangent to $B_\beta \cap U_F$
for all $\beta$ and to the fiber $F$ over $s$.
\end{prop}
\proof
That $\zeta_f$ is tangent to $B_\beta \cap U_F$
for all $\beta$ is the previous proposition.

To prove the second statement, note that
for every generic $b \in (F_\alpha)_{reg} \cap (B_\alpha)_{reg}$ and every $w \in T_b B_\alpha$, we have
$$\omega_b (\zeta_f (b), w) = d f_b (w) = d g_s (d \pi_b (w))=0$$
because $d \pi_b (w)$ is tangent to $S$ at $s$.
Hence $\zeta_f (b)$ belongs to the radical of $\omega|_B$.
By Theorem \ref{omega=0} (2), $\zeta_f (b)$ is tangent to $F$.
\endproof

Heuristically,
the  map $\pi: U_F \rightarrow {\Bbb C}^N$
can be interpreted as a {\it generalized}
moment map in the following sense.

There are many holomorphic
Hamiltonian flows around the fiber $F$.
Consider the coordinate functions $\{z_1, \ldots, z_N\}$.
They generate the Hamiltonian vector fields $\{\zeta_1, \ldots, \zeta_N\}$ which acts
(infinitesimally)
around a neighborhood of the fiber $F$ over $s$. We can think of
 $${\frak g} := \Span \{\zeta_1, \ldots, \zeta_N\}$$ as a `` Lie algebra''.
Its dual space ${\frak g}^* =\Span \{d z_1, \ldots,d  z_N \}$ is naturally
identified with ${\Bbb C}^N \cong m_{s, Z}/m_{s, Z}^2$.
Then the map
$$\pi: U_F \rightarrow {\Bbb C}^N \cong m_{s, Z}/m_{s, Z}^2 \cong {\frak g}^*$$
checks the conditions
$$d z_i = \iota_{\zeta_i} \omega, 1 \le i \le N$$
which characterize the most important feature of  a holomorphic moment map.
The presence of a genuine group action does not seem essential. What is important is that
many consequences of a moment map stems  only from the equations $d z_i = \iota_{\zeta_i} \omega$,
 $1 \le i \le N$.

 ${\frak g}$ has a `` Lie subalgebra''
$${\frak g}_s := \{\sum_i \lambda_i \zeta_i : \sum_i \lambda_i (d z_i)_s = 0 \; \hbox{along
$S$} \; ,
\lambda_i \in {\Bbb C} \}$$
$$= \{ g \in m_{s, Z}/m_{s, Z}^2 : (dg)_s = 0 \; \hbox{along $S$}. \;\}$$
This can be interpreted as the isotropy subalgebra at the point $s$.
Although ${\frak g}$ may not really be a Lie algebra, we suspect that
${\frak g}_s$ is indeed a Lie algebra.

We expect that the above Hamiltonian flows should play some important role in
birational geometry of complex (not necessarily projective) symplectic geometry. For example,
it is natural to speculate that
every irreducible component of a fiber of the small contraction $\pi: X \rightarrow Z$
is an almost homogeneous space.
Here, an arbitrary variety $V$ is said to be an almost homogeneous space of
a complex Lie group $G$ if $G$ acts on $V$ and has a dense open orbit.

Finally, we conjecture that

\begin{conj}
\label{codim2met} Let $\Phi: X ---> X'$ be a birational transformation between two smooth
projective holomorphic symplectic varieties of dimension $2n$.
Then, after removing subvarieties of codimension greater than 2,
$X$ and $X'$ are related by a sequence of Mukai's elementary transformations.
\end{conj}

An approach may go as follows.
 Let $E$ be the exceptional locus of $X$
and $E = \bigcup_i B_i$ be the union of its irreducible component.
We can arrange to have some components $\{B_j : j \in J\}$ contracted
 through log-extremal contraction.
Assume that $B$ is an arbitrary  one of the contracted components.
If $\codim B   >2$, we leave it alone. If $\codim B   =2$,
Theorem \ref{p2} indicates that
after throwing away a subvariety of codimension greater than 2 from $B$, the resulting
fibring
$B^0 \rightarrow S^0$ is a ${\Bbb P}^2$-bundle.

Now keep only these ${\Bbb P}^2$-bundles and remove all other subvarieties
in the contracted components  $\cup \{B_j : j \in J\}$
By Proposition \ref{disjoint},
all the above $\{B_j^0 : j \in J\}$ are disjoint.

The Mukai transformations can be performed along $\{B_j : j \in J\}$ (in any order).
Call the resulting symplectic variety $X^1$. One may repeat the above
to the birational map $X^1 ---> X'$.
After a finite step, we will arrive $X'$ from $X$ provided the necessary surgeries
are properly applied.

\vskip .4cm

\vskip .4cm

\bibliographystyle{amsplain}

\begin{thebibliography}{10}
\bibitem{BHL} D. Burns, Y. Hu, T. Luo,
{\em Symplectic Varieties and Birational Transformations in dimension 4,}
Contemporary Mathematics  of AMS (to appear).
\bibitem{Miyaoka} K. Cho, Y. Miyaoka,
{\em A characterization of Projective Spaces in Terms of the Minimum
Degrees of Rational Curves, }
Preprint (June 1997), RIMS, Kyoto University.
\bibitem{Miyaoka3} K. Cho, Y. Miyaoka, N. Shepherd-Barron,
{\em Characterizations of projective space and applications to
complex symplectic manifolds,} Adv. Studies in Pure Math. {\bf
35}, 2002, Higher dimensional birational geometry, pp. 1-88
\bibitem{fujiki} A. Fujiki,
{\em On primitively symplectic compact K\"ahler $V$-manifolds of dimension four.}
Classification of algebraic and analytic manifolds (Katata, 1982), 71--250,
Progr. Math., 39, Birkh\"auser Boston, Boston, MA, 1983.
\bibitem{huy1} D. Huybreschts,
{\em Birational symplectic manifolds and their deformations,}
J. Diff. Geom.
\bibitem{huy2} D. Huybreschts,
{\em Compact hyperK\"ahler manifolds: basic results,}
Invent. Math.
\bibitem{keladin} D. Kaledin,
{\em Symplectic resolutions: deformations and birational maps,}
math.AG/0012008.
\bibitem{kebekus} S. Kebekus,
{\em Characterizing the projective space after Cho, Miyaoka and
         Shepherd-Barron,} math.AG/0107069
\bibitem{Ko} J. Koll\'ar,
{\em Rational curves on algebraic varieties,}
Springer-Verlag, 1994.
\bibitem{km} J. Koll\'ar and S. Mori,
{\em Birational geometry of algebraic varieties,}
Cambridge University Press, 1998.
\bibitem{markman} E. Markman,
{\em Brill-Noether duality for moduli spaces of sheaves on K3 surfaces,}
 Journal of Algebraic Geometry 10 (2001), no. 4, 623-694.
\bibitem{Namikawa} Y. Namikawa,
{\em Deformation theory of singular symplectic n-folds,}
Math. Ann. {\bf 319} (2001), 597-623 (2001)).
\bibitem{SYZ} A. Strominger, S.-T. Yau, and E. Zaslow,
{\em Mirror Symmetry is T-Duality,}
Nucl. Phys. {\bf B479} (1996), 243-259.
\bibitem{wierzba} J. Wierzba,
{\em Contractions of Symplectic Varieties,} math.AG/9910130.
\bibitem{WW} J. Wierzba and J. Wisniewski,
{\em Small contractions of symplectic 4-folds}. math.AG/0201028
\end{thebibliography}
\makeatletter \renewcommand{\@biblabel}[1]{\hfill#1.}\makeatother

\end{document}